\documentclass{article}
\usepackage{amsmath,amssymb}
\usepackage{times}
\usepackage{graphicx}
\usepackage{color}
\usepackage{mathtools}

\def\spcc{{ \mbox{\hspace{.1in}}}}

\newcommand{\ld}{\lambda}
\newcommand{\frd}[2]{{\displaystyle \frac{#1}{#2}}}
\newcommand{\pdf}[2]{{\displaystyle \frac{\partial #1}{\partial #2}}}

\newcommand{\be}{\begin{eqnarray}}
\newcommand{\ee}{\end{eqnarray}}
\newcommand{\ben}{\begin{eqnarray*}}
\newcommand{\een}{\end{eqnarray*}}
\newcommand{\nn}{\nonumber}

\newcommand{\ta}{\mbox{\boldmath $a$}}

\newcommand{\tg}{\mbox{\boldmath $g$}}

\newcommand{\tb}{\mbox{\boldmath $b$}}
\newcommand{\tu}{\mbox{\boldmath $u$}}

\newcommand{\tv}{\mbox{\boldmath $v$}}

\newcommand{\tm}{\mbox{\boldmath $m$}}
\newcommand{\tn}{\mbox{\boldmath $n$}}
\newcommand{\tA}{\mbox{\boldmath $A$}}
\newcommand{\tB}{\mbox{\boldmath $B$}}
\newcommand{\tC}{\mbox{\boldmath $C$}}
\newcommand{\tD}{\mbox{\boldmath $D$}}
\newcommand{\tE}{\mbox{\boldmath $E$}}
\newcommand{\tF}{\mbox{\boldmath $F$}}
\newcommand{\tH}{\mbox{\boldmath $H$}}
\newcommand{\tI}{\mbox{\boldmath $I$}}
\newcommand{\tL}{\mbox{\boldmath $L$}}

\newcommand{\tQ}{\mbox{\boldmath $Q$}}
\newcommand{\tS}{\mbox{\boldmath $S$}}
\newcommand{\tT}{\mbox{\boldmath $T$}}

\newcommand{\tW}{\mbox{\boldmath $W$}}
\newcommand{\tU}{\mbox{\boldmath $U$}}

\newcommand{\tV}{\mbox{\boldmath $V$}}
\newcommand{\tG}{\mbox{\boldmath $G$}}

\newcommand{\bl}{\bullet}
\newcommand{\dis}{\displaystyle}

\newcommand{\tr}{\mbox{tr}}
\newcommand{\ru}[2]{ {#1}^{(#2)}}
\newcommand{\rp}[2]{ {#1}_{(#2)}}

\newcommand{\ot}{\otimes}

\def\rr#1{(\ref{#1})} 

\def \tr{\mbox{tr\hskip 1pt}}

\newtheorem{thrmc}{Theorem}

\def\theequation{\arabic{equation}}
\addtocounter{figure}{0}

\title{On the Smallest Number of Functions Representing Isotropic Functions of Scalars, Vectors and Tensors}
\author{M.H.B.M. Shariff\\[0.2cm]
\small Department of Applied Mathematics and Science\\
\small Khalifa University of Science, Technology and Research, UAE.}

\date{}

\begin{document}
\maketitle

\begin{abstract}
In this paper, we address the open problem (stated in Pennisi and Trovato, 1987. Int. J. Engng Sci., 25(8), 1059-1065) associated with the irreducibility of representations for isotropic functions. In particular, we prove that for isotropic functions
that depend on $P$ vectors, $N$ symmetric tensors and $M$ non-symmetric tensors
(a)  the number of irreducible invariants for a scalar-valued isotropic function is $3P+9M+6N-3$ (b)
the number of irreducible vectors for a vector-valued isotropic function is $3$ and (c)
the number of irreducible tensors for a tensor-valued isotropic function is at {\it most} $9$. 
The irreducible numbers in given (a), (b) and (c) are much lower than those obtained in the literature. This significant reduction in the number of
irreducible scalar/vector/tensor-valued functions have the potential to substantially simplify modelling complexity.
\end{abstract}

\section{Introduction}
Mathematical modelling of physical conditions often requires representations for isotropic functions \cite{pipkin58,rivlin55}. In view of this much has been published
on this subject (see, for example reference \cite{penninsi87}, and references therein). However, the derived number of isotropic functions in an irreducible basis (see definition of an irreducible basis in \cite{spencer71})  is still an open problem as stated by
Pennisi and Trovato \cite{penninsi87}, where they state that: \\
"{\it  Among all  irreducible complete representations previously published in the literature (2.1)-(2.4) is that with fewer elements;
but it is still an open problem to find, among all {\bf possible} irreducible complete representations, that (if it exists) with fewer elements}".

In this paper, we address this open problem and prove that only a few elements are required  in irreducible bases.
The proofs given here are simple (compared to the proofs given in the literature) and they are based on a spectral approach 
associated with the author's work \cite{shariff13,sharbusta15,shariff17,shariff21a}. This substantial reduction in numbers  of elements in irreducible bases could radically reduce modelling complexity.

\section{Preliminaries}
Let $V$  be a $3$-dimensional vector space.
We define  $Lin$ to be the space of all linear
transformations (second-order tensors) on $V$ with the inner product $\tA:\tB=tr(\tA\tB^T)$, where $\tA, \tB \in Lin$ and $\tB^T$ is the transpose
of $\tB$. We define
\be\label{pra} Sym=\{\tA\in Lin | \tA=\tA^T\} \, , \spcc Orth= \{ \tQ \in Lin | \tQ = \tQ^{-T} \} \, . \ee
The vectors considered here belong to the $3$-dimensional Euclidean space $\mathbb{E}^3$, i.e., the vector space $V$ furnished by the
scalar product $\ta\cdot\tb$, where $\ta,\tb \in V$.

The summation convention is not used here  and, all subscripts $i, j$ and $k$ take the values $1,2,3$ unless stated otherwise.

\section{Symmetric Tensors and Vectors}

\subsection{Scalar}\label{sec-scalar-1}
The scalar function $W(\tA_r,\ta_s)$, $(r=1,2,\ldots , N ; s=1,2,\ldots , P)$, where $\tA_r \in Sym $ and $\ta_s  \in \mathbb{E}^3$ are, respectively,
symmetric tensors and vectors, is said to be scalar-valued isotropic function if
\be W(\tA_r,\ta_s) = W(\tQ\tA_r\tQ^T,\tQ\ta_s) \, \ee
for all rotation tensor $\tQ \in Orth$.  Boehler \cite{boehler77} has shown that every scalar-valued
isotropic function can be written as a function of invariants given in the following list:
\[ \ta_\alpha\cdot\ta_\alpha \, , \spcc \ta_\alpha\cdot\ta_\beta \, , \]
\[ \tr \tA_i\, , \spcc \tr \tA_i^2 \, , \spcc \tr \tA_i^3 \spcc \tr \tA_i^2\tA_j \, , \spcc \tr \tA_i\tA_j^2 \, , \spcc
\tr \tA_i^2\tA_j^2 \, , \spcc \tr \tA_i\tA_j\tA_k \, , \]
\[ \ta_\alpha\cdot\tA_i\ta_\alpha \, , \spcc \ta_\alpha\cdot\tA_i^2\ta_\alpha \, , \spcc 
\ta_\alpha\cdot\tA_i\tA_j\ta_\alpha \, , \]
\be\label{sm-1}\ta_\alpha\cdot\tA_i\ta_\beta \, , \spcc  \ta_\alpha\cdot\tA_i^2\ta_\beta \, , \spcc 
\ta_\alpha\cdot(\tA_i\tA_j - \tA_j\tA_i) \ta_\beta \, , \ee
$i,j,k=1,2,\ldots , N$ with $i<j<k$ and $\alpha,\beta= 1,2,\ldots , P$ with $\alpha < \beta$. However, Shariff \cite{shariff21a} has shown
that, for unit vectors $\tv_\alpha$, only $2P+6N-3$ of the invariants in \rr{sm-1} are independent and that the number of invariants in
the irreducible functional basis  is at most $2P+6N-3$; far lower than the number of invariants given in \rr{sm-1}. In the case when $\tv_\alpha$ are not
unit vectors it can be easily shown that only $3P+6N-3$ of the invariants in \rr{sm-1} are independent. Below, for the sake of easy reading, 
we prove (similar to the work of Shariff \cite{shariff21a}) that every scalar-valued
isotropic function can be written as a function of at most $3P+6N-3$ number of invariants. This significant reduction in number  of scalar invariants
(when compared to the list in \rr{sm-1}) could greatly assist in reducing modelling complexity (see for example references 
\cite{shariff13a,shariff14,shariff16,shariff17b,shariff20a,shariff20b,shariff21b,shariff22a,shariff22b,shariff22c})\\
{\it Proof}\\
For $N \ge 1$. Let express (say)
\be\label{scalar-1} \tA_1 = \sum_{i=1}^3 \ld_i \tv_i\ot\tv_i \, , \ee
where $\ld_i$ and $\tv_i$ are eigenvalues and (unit) eigenvectors of $\tA_1$, respectively and $\ot$ represents a dyadic product. Using $\{\tv_1,\tv_2,\tv_3 \}$ as a basis,
we can express
\be \tA_r = \sum_{i,j=1}^3 \ru{A}{r}_{ij} \tv_i\ot\tv_j \, , \spcc \ta_s = \sum_{i=1}^3 \ru{a}{s}_i \tv_i \, , \spcc
r=2,3, \ldots N \, , \spcc s=1,2,\ldots , P \, . \ee
It is clear that the components $\ru{A}{r}_{ij}$ and $\ru{a}{s}_i$ are invariants, since
\be \ru{A}{r}_{ij} = \tv_i\cdot\tA_r\tv_j = \tQ\tv_i\cdot\tQ\tA_r\tQ^T\tQ\tv_j \, , \spcc
\ru{a}{r}_i = \ta_r\cdot\tv_i = \tQ\ta_r \cdot \tQ\tv_i \, . \ee
Since,
\be\label{ire-1} \ld_i \, , \spcc \ru{A}{r}_{ij} \, , \spcc \ru{a}{s}_i \, , \spcc r \ge 2 \, , \spcc i,j=1,2,3 \,    \ee
are "component" invariants, we can express
\be W(\tA_r,\ta_s) = W(\tQ\tA_r\tQ^T,\tQ\ta_s) = \hat{W} (\ld_i,\ru{A}{r}_{ij},\ru{a}{s}_i ) \, , \spcc r \ge 2  \spcc i,j=1,2,3 \, .  \ee
All invariant functions in \rr{sm-1} can be explicitly expressed in terms of the  spectral invariants given below;
for example, we can express the function
\be \ta_\alpha\cdot\tA_i^2\ta_\beta = \sum_{p,q,m=1}^3 \ru{a}{\alpha}_p\ru{A}{i}_{pq}\ru{A}{i}_{qm}\ru{a}{\beta}_m \, , \spcc i \ne 1 \ee
Hence, the set of invariants in \rr{ire-1} is a complete representation for the scalar-valued isotropic function and 
since the terms in \rr{ire-1} are  independent (invariant) components, the set is 
irreducible, i.e., incapable of being reduced.
Hence, every scalar-valued isotropic function can be written as a function of at most $3P+6N-3$ number of invariants,
far less than the number of invariants given in \rr{sm-1}.
The  spectral invariants in \rr{ire-1} have been used in continuum modelling \cite{shariff13a,shariff14,shariff16,shariff17b,shariff20a,shariff20b,shariff21b,shariff22a,shariff22b,shariff22c} 
and spectral derivatives, associated with these spectral invariants, are given in \cite{shariff17a,shariff20}.

Since all of Boehler's  invariants \rr{sm-1} can be explicitly expressed in terms of the spectral invariants \rr{ire-1}, this further validate
our claim that the irreducible basis contains only $6N +3P-3$ invariants.

{\bf Word of caution:} The function 
\be \hat{W}(\ld_i,\ru{A}{r}_{ij},\ru{a}{s}_i )   \ee
must satisfy the $P$-property given in \cite{shariff16} and (for the benefit of the readers) in Appendix A. In this paper, we call
a scalar-valued isotropic function that satisfies the $P$-property, a $P$-scalar-valued isotropic function. In general, the invariants
appearing (as they are) in \rr{ire-1} are not $P$-scalar-valued isotropic functions.

In the case when $N=0$, we have $W$ depends on $\ta_s$ only. In this case,
we select the vector $\ta_1$ (say) and spectrally express
\be\label{vec-a1}
\ta_1 \ot\ta_1 = \ld\tv_1\ot\tv_1 + 0 \tv_2\ot\tv_2 + 0 \tv_3\ot\tv_3 \, , \spcc \ld = \ta_1\cdot\ta_1 \, , \spcc \tv_1 = \frd{\ta_1}{\sqrt{\ld}} \, \ee
and, $\tv_2$ and $\tv_3$ are any two (non-unique) orthonormal vectors that are perpendicular to $\ta$. Hence, for $N=0$, we have $3P-2$ 
irreducible invariants, i.e.,
\be \ld \, , \spcc \ru{a}{s}_i \, , \spcc s=2,3,\ldots , P \, , \spcc i=1,2,3 \, .  \ee
In the case where all of the vectors $\ta_s$ are unit vectors, we have only $2P-2$ irreducible spectral invariants.\\
{\bf Example 1:} Consider the strain energy function $W$ of a transversely isotropic elastic solid. We then have, 
\be\label{tr-1} W(\tU,\ta\ot\ta) = \tilde{W}(\tU,\ta) = \hat{W}(\ld_i, a_i) \, , \spcc a_i =\tv_i\cdot\ta \, , \ee
where $\ta_1=\ta$ is the preferred direction unit vector, $\tA_1=\tU$ is the right stretch tensor and
\be\label{tr-2} \sum_{i=1}^3 a_i^2 =1 \, . \ee
It is clear from \rr{tr-1} and \rr{tr-2}, and if we consider the positive and negative
values of $a_i$ as distinct single-valued functions then we can conclude that the number of invariants
in the irreducible functional basis is $5$.\\
{\bf Example 2:} If we consider in Example 1, $\tA_1=\ta\ot\ta$ and $\tA_2=\tU$, we have 
\be \ld_1 = 1 \, , \spcc \ld_2=\ld_3=0 \, , \spcc \tv_1=\ta \, ,  \ee
$\tv_2$ and $\tv_3$ are any two (non-unique) orthonormal vectors that are perpendicular to $\ta$ and we then have
\be\label{tr-3}  W(\ta\ot\ta,\tU) = \hat{W}(U_{ij}) \, , \spcc U_{ij} = \tv_i\cdot\tU \tv_j \, . \ee
We note that there are $6$ (instead of $5$) spectral invariants in \rr{tr-3}. However, since $\hat{W}$ must satisfy the $P$-property,
we can express $\hat{W}$ in terms of $5$ independent invariants, that satisfy the $P$-property. For example, we can express $\hat{W}$ in terms of the 5 independent invariants
\be I_1 = \sum_{i=1}^3 U_{ii} \, , \spcc I_2 = \sum_{i,j=1}^3 U_{ij}U_{ji} \, , \spcc I_3 = \sum_{i,j,k=1}^3 U_{ij}U_{jk}U_{ki} \, , \spcc
I_4= U_{11} \, , \spcc  I_5 = \sum_{i=1}^3 U_{1i}U_{i1} \, . \ee

\subsection{Vector}\label{sec-vector-1}
The vector function $\tg(\tA_r,\ta_s)$ is said to be vector-valued isotropic function if
\be\label{vect-1} \tQ\tg(\tA_r,\ta_s) = \tg(\tQ\tA_r\tQ^T,\tQ\ta_s) \, \ee
for all rotation tensor $\tQ$.

Smith \cite{smith71} has shown that every vector-valued isotropic function can be written as a linear combination of the following
vectors
\be\label{sm-2} \ta_m \, , \spcc \tA_i\ta_m \, , \spcc \tA_i^2\ta_m \, , \spcc (\tA_i\tA_j - \tA_j\tA_i)\ta_m \, , 
\spcc i,j=1,2,\ldots , N \, : i<j \, , \spcc m=1,2, \ldots , P \, . \ee
It is understood that the coefficients in these linear combinations are $P$-scalar-valued isotropic functions.

Smith \cite{smith71}  and Pennisi and Trovato  \cite{penninsi87} claimed that the set of vectors in \rr{sm-2} is irreducible; we claim that the irreducible set contains only  three linearly independent vectors. Below, we show via a theorem  that every vector-valued isotropic function can be written as a linear combination of at most three linearly independent spectral vectors. \\
\begin{thrmc} $\tg$ is an isotropic tensor function if and only if it has the representation
\be\label{apc-2} \tg(\tA_r,\ta_s) = \sum_{i=1}^3 g_i \tv_i \, , \ee
where $\tv_i$ is an eigenvector of $\tA_1$ and $g_i$ are isotropic invariants of the  set 
\be\label{apc-2a} S = \{\tA_1,\tA_2, \ldots \tA_N, \ta_1,\ta_2 , \ldots, \ta_P \} \, . \ee
\end{thrmc}
{\it Proof:} \\
(a) If \rr{apc-2} holds $\tg$ is clearly a vector-valued isotropic function, since the coefficients $g_i$ are isotropic invariants of the set
$S$ \rr{apc-2a}. \\
(b) For $N \ge 1$ and $P\ge 0$. Let $\tv_i$ be unit eigenvectors of the symmetric tensor $\tA_1$ (see \rr{scalar-1}).
Hence we can write
\be\label{vect-2} \tg(\tA_r,\ta_s) = \sum_{i=1}^3 [\tg(\tA_r,\ta_s) \cdot\tv_i]\tv_i  \,  ,  \spcc r=1,2,\ldots, N, \spcc
s=1,2,\ldots , P \ee
and
\be\label{vect-3} \tg(\tQ\tA_r\tQ^T,\tQ\ta_s) = \sum_{i=1}^3 [\tg(\tQ\tA_r\tQ^T,\tQ\ta_s)\cdot\tQ\tv_i] \tQ\tv_i  \,  . \ee
Let scalar function
\be g_i(\tA_r,\ta_s) = \tg(\tA_r,\ta_s) \cdot\tv_i \, . \ee
We then have
\be g_i(\tQ\tA_r\tQ^T,\tQ\ta_s) = \tg(\tQ\tA_r\tQ^T,\tQ\ta_s)\cdot\tQ\tv_i  \, . \ee
In view of \rr{vect-1}, \rr{vect-2} and \rr{vect-3}, and since $\tQ$ is arbitrary, we must have
\be g_i(\tA_r,\ta_s) = g_i(\tQ\tA_r\tQ^T,\tQ\ta_s) \, , \ee
which implies that $g_i$ are functions of isotropic invariants of the vector and tensor set $S$ given in \rr{apc-2a}.
Note that, in view of the $P$-property,  the functions $g_i$ must also be $P$-scalar-valued isotropic functions.

In the case when $N=0$, we consider the vectors $\tv_i$ obtained similar to \rr{vec-a1} and express
\be \tg(\ta_r) = \sum_{i=1}^3 g_i \tv_i \, ,  \spcc g_i = \tg\cdot\tv_i \, .  \ee
All Smith's vectors given in \rr{sm-2} can be expressed in terms of the unit vectors $\tv_1,\tv_2$ and $\tv_3$. For example the vector
\be \tA_i\ta_m = \sum_{r=1}^3 (\sum_{s=1}^3\ru{A}{i}_{rs}\ru{a}{m}_s)\tv_r \, . \ee
Hence, when a vector-valued function is expressed in terms of a linear combinations of Smith's functions given in \rr{sm-2}, it can then be
expressed in terms of a linear combination of the symmetric spectral vectors $\tv_1,\tv_2$ and $\tv_3$; this further validates our claim
that the irreducible basis contains only three vectors.
\subsection{Symmetric Tensor}\label{sec-tensor-1}
The symmetric tensor function $\tG(\tA_r,\ta_s)$, is said to be tensor-valued isotropic function if
\be\label{tensor-1} \tQ\tG(\tA_r,\ta_s)\tQ^T = \tG(\tQ\tA_r\tQ^T,\tQ\ta_s) \, \ee
for all rotation tensor $\tQ$.

Smith \cite{smith71} has shown that every symmetric tensor-valued isotropic function can be written as a linear combination of the following
symmetric tensors
\[ \tI \, , \spcc \tA_i \, , \spcc \tA^2_i \, , \spcc \tA_i\tA_j+ \tA_j\tA_i \, , \spcc
\tA_i^2\tA_j+ \tA_j\tA_i^2 \, , \spcc \tA_i\tA_j^2+ \tA_j^2\tA_i \]
\[ \ta_m\ot\ta_m \, , \spcc \ta_m\ot\ta_n +\ta_n\ot\ta_m \, , \spcc
\ta_m\ot\tA_i\ta_m + \tA_i\ta_m\ot\ta_m \, , \spcc  \ta_m\ot\tA_i^2\ta_m + \tA_i^2\ta_m\ot\ta_m \, , \]
\be\label{sm-3}
\tA_i(\ta_m\ot\ta_n - \ta_n\ot\ta_m) - (\ta_m\ot\ta_n - \ta_n\ot\ta_m)\tA_i \, , \ee
where $(i,j=1,2, \ldots , N; i<j)$, $(p,q = 1,2, \ldots , M ; p < q)$, $(m,n =1,2, \ldots , P; m< n)$ and $\tI$ is the
identity tensor.
Smith \cite{smith71}  and Pennisi and Trovato  \cite{penninsi87} claimed that the set of symmetric tensors in \rr{sm-3} is irreducible; we,
however, claim
via Theorem \ref{thrm-2} below, that the irreducible set contains only  six linearly independent symmetric tensors.
\begin{thrmc}\label{thrm-2} $\tG$ is an isotropic tensor function if and only if it has the representation
\be\label{apc-2ten} \tG(\tA_r,\ta_s) = \sum_{i,j=1} t_{ij} \tv_i\ot \tv_j \, , \ee
where $\tv_i$ is an eigenvector of $\tA_1$ and $t_{ij}$ are functions of $P$-scalar-valued isotropic functions of the vector and tensor set 
given in \rr{apc-2a}.
\end{thrmc}

{\it Proof} \\
(a) If \rr{apc-2ten} holds, since $t_{ij}$ are scalar invariants of the set $S$, then $\tG$ is clearly and isotropic tensor function. \\
(b) Using the basis $\{ \tv_1,\tv_2,\tv_3 \}$ obtained from \rr{scalar-1} , we can express
\be \tG(\tA_r,\ta_s) = \sum_{i,j=1} t_{ij} \tv_i\ot\tv_j \, , \ee
where
\be t_{ij} = g_{ij}(\tA_r,\ta_s) = \tv_i\cdot\tG(\tA_r,\ta_s)\tv_j \,  . \ee
Similarly, we can express 
\be \tG(\tQ\tA_r\tQ^T, \tQ\ta_s)
=\sum_{i,j=1} \bar{t}_{ij} \tQ\tv_i\ot\tQ\tv_j  \, , \ee
where
\be &&\bar{t}_{ij} = \tQ\tv_i\cdot\tG(\tQ\tA_r\tQ^T, \tQ\ta_s)\tQ\tv_j \nn \\
&& =g_{ij}(\tQ\tA_r\tQ^T,\tQ\ta_s) \, . \ee
If \rr{tensor-1} holds then
\be \sum_{i,j=1} \bar{t}_{ij} \tQ\tv_i\ot\tQ\tv_j = \sum_{i,j=1} t_{ij} \tQ\tv_i\ot\tQ\tv_j  \, . \ee
Since $\tQ$ is arbitrary, we have
\be g_{ij}(\tA_r,\ta_s) = g_{ij}(\tQ\tA_r\tQ^T, \tQ\ta_s) \, \ee
which implies that the functions $t_{ij}=g_{ij}$ must depend on $P$-scalar-valued isotropic functions of $S$. Since
$g_{ij}=g_{ji}$, all tensor-valued isotropic functions can be written as  a linear combination of only six symmetric tensors 
\be\label{gv-sm3}  \tv_i\ot\tv_i \,  \spcc (i=1,2,3) \, , \spcc  \tv_i\ot\tv_j+\tv_j\ot\tv_i \, \spcc( i=1,2; j=2,3, i<j) \, . \ee
Hence, we can express
\be\label{sym-ga} \tG(\tA_r, \ta_s) =\sum_{i=1}^3 g_{ii} \tv_i\ot\tv_i + \sum_{i<j} g_{ij}(\tv_i\ot\tv_j+\tv_j\ot\tv_i) \, \spcc (i=1,2;j=2,3) \, . \ee
All symmetric tensors in \rr{sm-3} generated by Smith \cite{smith71} can be expressed in terms of the six symmetric tensors given in 
\rr{gv-sm3}, for example, the symmetric tensor
\be \tA_i\tA_j+ \tA_j\tA_i = \sum_{p=1}^3 g_{pp} \tv_p\ot\tv_p + \sum_{p<q} g_{pq}(\tv_p\ot\tv_q+\tv_q\ot\tv_p) \, \spcc (p=1,2;q=2,3) \, , \ee
where
\be g_{pp} = 2\sum_{m=1}^3 \ru{A}{i}_{pm}\ru{A}{j}_{mp} \, , \spcc g_{pq} = \sum_{m=1}^3 (\ru{A}{i}_{pm}\ru{A}{j}_{mq}+
\ru{A}{i}_{qm}\ru{A}{j}_{mp} ) \, . \ee
Hence, when a tensor-valued function is expressed in terms of a linear combinations of Smith's functions given in \rr{sm-3}, it can then be
expressed in terms of a linear combination of the six symmetric spectral tensors given in \rr{gv-sm3}; this further validates  our claim
that the irreducible basis contains only six symmetric tensors.

The above theorem proves that the irreducible set contains only  six linearly independent symmetric tensors. This drastically
reduce the complexity in physical modelling. For example, Merodio and Rajagopal \cite{mer07} modelled viscoelastic solids, where the Cauchy stress
$\tT$ depends on  $\tA_1=\tB$ (left Cauchy-Green
stretch tensor), $\tA_2=\tD$ (the symmetric part of the velocity gradient), $\ta_1=\tm$ and $\ta_2=\tn$ (preferred directions). Using 
Smith tensors \rr{sm-3}, the Cauchy stress $\tT$ is described using $36$ tensors obtained from \rr{sm-3} and,  due to this large number of
$36$ tensors and $37$ scalar invariants, the model is complicated; there is a dire need to simplify the model. Sometimes this is done by
omission of invariants and tensors. However,
the discrimination in selection of invariants and tensors  is often debated, and neglecting the influence of some invariants and tensors may
result in an incomplete representation of the full range of mechanical response subjected to a continuum. However, using the results obtained here,
modelling viscoelastic solids is greatly simplified, we only require $15$ scalar invariants and $6$ symmetric tensors to fully describe
the Cauchy stress $\tT$.

{\bf Remark:}\\
{\it Since both the scalars $g_i$ and $g_{ij}$ are, respectively, vector and tensor components, the vector $\tg$ and 
tensor $\tG$ are uniquely expressed in terms of the basis $\{\tv_1,\tv_2,\tv_3 \}$ even though two or three of the vectors
$\tv_1$, $\tv_2$ and $\tv_3$ may not be unique due to coalescence of eigenvalues.}

The theorem below has been proven in the literature (see for example references Itskov \cite{itskov2013} and Ogden \cite{ogden84}), however, for 
the benefit of the readers we prove it again here.
\begin{thrmc}
If $\tG(\tV)$ is an isotropic tensor function then $\tG(\tV)$ is coaxial with $\tV$ and hence
\be \tV\tG = \tG\tV \, . \ee
\end{thrmc}
{\it Proof} \\
Let $\tv_1$ be an eigenvector of $\tV$ and choose
\be \tQ = 2\tv_1\ot\tv_1 - \tI = \tQ^T \, . \ee
In view of $\tV\tv_1 = \ld_1 \tv_1$, we have
\be  \tQ\tV = \tV\tQ \rightarrow \tQ\tV\tQ^T = \tV \, . \ee
From \rr{tensor-1} we get
\be \tQ\tG(\tV)\tQ^T = \tG(\tV) \, . \ee
Hence
\be \tG(\tV)\tv_1 = \left( \sum_{i,j=1} t_{ij} \tQ\tv_i\ot\tv_j\tQ^T \right) \tv_1 \, . \ee
Note that $\tQ^T\tv_1 = \tv_1$ and hence we have
\be \tG(\tV)\tv_1 =  \sum_{i=1} t_{i1} \tQ\tv_i = \sum_{i=1} t_{i1} (2\tv_1\ot\tv_1 - \tI)\tv_i  = 2t_{11}\tv_1 - \tG(\tV)\tv_1 \, .  \ee
Hence
\be \tG(\tV)\tv_1 = t_{11}\tv_1 \, , \ee
which implies that $\tv_1$ is an eigenvector of $\tG(\tV)$ and $t_{11}$ is an eigenvalue of $\tG$. In a similar fashion,
choosing $\tQ = 2\tv_r\ot\tv_r - \tI = \tQ^T$, $r=2,3$, we can easily derive that
\be \tG(\tV) = \sum_{i=1} t_{ii} \tv_i\ot\tv_i \, ,\ee
and the theorem is proved. 

Below is a theorem, which we believe is not found in the literature.
\begin{thrmc}
Let $\ld_i$ be the eigenvalues of $\tV$ and let
\be \tG(\tV) = \sum_{i=1}^3 t_i (\ld_1,\ld_2,\ld_3) \tv_i\ot\tv_i \,   , \ee
be a symmetric isotropic tensor function, where $\tv_i$ is an eigenvector of $\tV$. \\
(a) If $\ld_i=\ld_j \ne \ld_k$, $(i \ne j \ne k \ne i)$, then
\be t_i=t_j \,  \ee
and we can uniquely express
\be \tG(\tV) = t_i\tI + (t_k-t_i)\tv_3\ot\tv_3 \, . \ee
(b)  If $\ld_1=\ld_2 =\ld_3$ then
\be t_1=t_2=t_3 \,  \ee
and we can uniquely express
\be \tG(\tV) = t_1\tI \, . \ee
\end{thrmc}
{\it Proof} \\
Consider the case $\ld_1=\ld_2 = \ld \ne \ld_3$. In view of this, $\tv_1$ and $\tv_2$ are not unique and have infinitely many values. 
In view of the relation
\be \tv_1\ot\tv_1 + \tv_2\ot\tv_2+\tv_3\ot\tv_3 = \tI \, , \ee
we can write
\be\label{GV-1} \tG(\tV) = t_1\tI + (t_2-t_1)\tv_2\ot\tv_2 + (t_3-t_1)\tv_3\ot\tv_3 \, . \ee
Since $\tv_2$ is not unique, we must have $t_1=t_2$ to give $\tG(\tV)$ a unique value. In a similar fashion, 
we can show for the cases $\ld_1=\ld_3$ and $\ld_2=\ld_3$. Hence, theorem (a) is proved.

In the case when $\ld_1=\ld_2 =\ld_3$, $\tv_3$ is also arbitrary, hence from \rr{GV-1} we must have $t_1=t_2=t_3$
and theorem (b) is proved.

We can see that in case when the classical invariants $I_1=\tr \tV$, $I_2=\tr \tV^2 $ and $I_3 = \tr \tV^3$ are used, we have \cite{ogden84}
\be\label{GV-2} \tG(\tV) = \phi_0 \tI + \phi_1 \tV + \phi_2\tV^2  = \sum_{i=1}^3 (\phi_0 + \phi_1\ld_i + \phi_2 \ld_i^2) \tv_i\ot\tv_i =
 \sum_{i=1}^3 t_i \tv_i\ot\tv_i \,  , \ee
\be\label{GV-3} t_i= \phi_0 + \phi_1\ld_i + \phi_2 \ld_i^2 \, , \ee
where $\phi_0,\phi_1$ and $\phi_2$ depend on  $P$-scalar-valued isotropic functions, $I_1$, $I_2$ and $I_3$. It is clear from \rr{GV-3} that $t_i=t_j$ when $\ld_i=\ld_j$.

\section{Isotropic Functions of Non-symmetric Tensors}
\subsection{Scalar}
The scalar function $W(\tH_t, \tA_r,\ta_s)$, $(r=1,2,\ldots , N ; t=1,2,\ldots M ; s=1,2,\ldots , P )$ is said to be a scalar-valued isotropic function if
\be W(\tH_t, \tA_r,\ta_s) = W(\tQ\tH_t\tQ^T, \tQ\tA_r\tQ^T,\tQ\ta_s) \, \ee
for all rotation tensor $\tQ \in Orth$, where $\tH_t \in Lin$ $(t=1,2,\ldots M)$ is a nonsymmetric second order tensor. 

In the case when $M,N,\ge 1$, we can easily proved, based on Section \ref{sec-scalar-1} that
\be W(\tH_t, \tA_r,\ta_s) = W(\tQ\tH_t\tQ^T, \tQ\tA_r\tQ^T,\tQ\ta_s) = \hat{W}(\ld_i, \ru{H}{t}_{ij},\ru{A}{r}_{ij}, \ru{a}{s}_i) \, ,
\spcc  r=2,3,\ldots N \, ,\ee
where the invariants  
\be\label{nsym-1a} \ld_i,\ru{A}{r}_{ij}, \ru{a}{s}_i \ee 
are given in \rr{ire-1} and the invariants
\be\label{nsym-1b} \ru{H}{t}_{ij} = \tv_i \cdot \tH_t \tv_j =\tQ\tv_i \cdot\tQ\tH_t\tQ^T \tQ\tv_j\, , \spcc i,j=1,2,3 \, . \ee
Since the above invariants are independent components, the irreducible basis consists of at most
$3P+9M+6N-3$ invariants. Note that Boehler \cite{boehler77}  consider the isotropic function
\be\label{nsym-2}  W(\tW_t, \tA_r,\ta_s)\, , \ee
where $\tW_t$ is a skew-symmetric tensor. He claimed that the irreducible set contains the "complicated" set of invariants
\[ \ta_\alpha\cdot\ta_\alpha \, , \spcc \ta_\alpha\cdot\ta_\beta \, , \spcc \tr \tA_i\, , \spcc \tr \tA_i^2 \, , \spcc \tr \tA_i^3 \, ,
\, , \spcc \tr \tA_i\tA_j \, , \spcc \tr \tA_i^2\tA_j \, , \spcc \tr \tA_i\tA_j^2 \, , \spcc \tr \tA_i^2\tA_j^2 \, ,\]
\[ \tr \tA_i\tA_j\tA_k \, , \spcc
 \tr \tW_p^2 \, , \spcc \tr \tW_p\tW_q \, , \spcc \tr \tW_p\tW_q\tW_r \, , \spcc
\ta_\alpha\cdot\tA_i \ta_\alpha \, , \spcc \ta_\alpha\cdot\tA_i^2 \ta_\alpha \, , \] 
\[ \ta_\alpha\cdot\tA_i\tA_j \ta_\alpha 
\, , \spcc \ta_\alpha\cdot\tA_i \tv_\beta \, , \spcc \ta_\alpha\cdot\tA_i^2 \ta_\beta \, , \]
\[ \ta_\alpha\cdot (\tA_i\tA_j-\tA_j\tA_i)\ta_\beta \, , \spcc \ta_\alpha\cdot\tW_p^2\ta_\alpha
\, , \spcc \ta_\alpha\cdot\tW_p\tW_q\ta_\alpha \, , \spcc \ta_\alpha\cdot\tW_p^2\tW_q\ta_\alpha \, ,\]
\[ \ta_\alpha\cdot\tW_p\tW_q^2\ta_\alpha \, , \spcc
 \ta_\alpha\cdot\tW_p\ta_\beta \, , \spcc  \ta_\alpha\cdot\tW_p^2 \ta_\beta \, , \]
\[ \ta_\alpha\cdot (\tW_p\tW_q-\tW_q\tW_p)\ta_\beta 
\, , \spcc \tr \tA_i\tW_p^2 \, , \spcc  \tr \tA_i^2\tW_p^2 \, , \spcc  \tr \tA_i^2\tW_p^2\tA_i\tW_p \, , \spcc 
\tr \tA_i\tW_p\tW_q \, , \spcc \tr \tA_i\tW_p\tW_q^2 \, , \]
\[  \tr \tA_i\tW_p^2\tW_q \, , \spcc
\tr \tA_i\tA_j\tW_p \, , \spcc \tr \tA_i\tW_p^2\tA_j\tW_p \, , \spcc \tr \tA_i\tA_j^2\tW_p 
\, , \spcc\tr \tA_i^2\tA_j\tW_p \, ,\]
\be\label{nysm-1}
\ta_\alpha\cdot\tA_i\tW_p \ta_\alpha \, , \spcc \ta_\alpha\cdot\tW_p\tA_i\tW_p^2 \ta_\alpha
\, , \spcc \ta_\alpha\cdot\tA_i^2\tW_p \ta_\alpha \, , \spcc 
 \ta_\alpha\cdot (\tA_i\tW_p-\tW_p\tA_i)\ta_\beta \, , \ee
where $i,j,k=1,2, \ldots , N$ with $i<j<k; p,q,r =1,2,\ldots , M$ with $p <q < r$ and
$\alpha,\beta =1,2,\ldots , P$  with $\alpha < \beta$. However, prove that irreducible set contains only
$3P+3M+6N-3$ invariants and they are:
\be\label{nsym-3} \ld_i\, , \spcc \ru{A}{r}_{ij}\, , \spcc  \ru{a}{s}_i \, , \spcc \ru{W}{t}_{kl} = \tv_k\cdot\tW_k\tv_l \, , \spcc i,j,k,l=1,2,3,\, ,
\spcc k<l \, , \spcc r \ge 2 \, . \ee
The invariants in \rr{nsym-3} are  obtained from \rr{nsym-1a} and \rr{nsym-1b}, by replacing $\tH_t$ with $\tW_t$ and
taking note that
\be \tv_i\cdot\tW_k\tv_i = 0 \, , \spcc \tv_i\cdot\tW_k\tv_j = -\tv_j\cdot\tW_k\tv_i \, , \spcc i\ne j \, , \spcc i,j=1,2,3 \, . \ee
In the case when $N=0$, we have
\be W(\tH_t,\ta_s) \, . \ee
In this case, we let  the orthonormal vectors $\tv_i$ to be the eigenvectors of the symmetric tensor $\tH_1\tH_1^T$ (or alternatively $\tH_1^T\tH_1$), i.e.,
\be \tH_1\tH_1^T = \sum_{i=1}^3 \ld_i \tv_i\ot\tv_i \, , \spcc \ld_i \ge 0 \, . \ee
The irreducible set contains at most $9M+3P$ invariants
\be\label{scalarb-1} \ru{H}{t}_{ij} \, , \spcc \ru{a}{s}_i \, , \spcc i,j=1,2,3 \, . \ee

\subsection{Vector}
For a vector-valued isotropic function, it can be easily prove that, following Section \ref{sec-vector-1},
\be \tg(\tH_t, \tA_r,\ta_s) = \sum_{i=1}^3 g_i \tv_i \, , \ee
where $g_i$ are functions of the invariants in \rr{nsym-1a} and \rr{nsym-1b} or \rr{scalarb-1}, as appropriate.
Hence, the irreducible basis for $\tg$ contain only the three vectors $\tv_i$. Note that for $\tH_t=\tW_t$, Smith \cite{smith71} claimed that
the irreducible basis for $\tg$ contain the vectors
\[ \ta_m \, , \spcc \tA_i\ta_m \, , \spcc \tA_i^2\ta_m \, , \spcc (\tA_i\tA_j - \tA_j\tA_i)\ta_m \, , \spcc
\tW_p\ta_m \, , \spcc \tW_p^2\ta_m \, , \]
\be\label{smvec-2}  (\tW_p\tW_q -\tW_q\tW_p)\ta_m \, , \spcc (\tA_i\tW_p - \tW_p\tA_i)\ta_m \, , \ee
where $i,j=1,2,\ldots , N \, ; i<j \, , \spcc p,q=1,2,\ldots , M; p<q \, ,  \spcc m=1,2, \ldots , P$.
This claim is incorrect since all the vectors in \rr{smvec-2} can be written terms of the vectors $\tv_1,\tv_2$ and $\tv_3$.

\subsection{Tensor}
Following the method in Section \ref{sec-tensor-1}, we can easily prove that the for any tensor in $Lin$, with $M,N\ge 1$,
\be \tH(\tH_t, \tA_r,\ta_s) = \sum_{i,j=1}^3 h_{ij} \tv_i\ot\tv_j \, , \ee
where, $\tv_i$ is an eigenvector of $\tA_1$, and in general $h_{ij}=\tv_i\cdot\tH\tv_j \ne h_{ji}$ are functions of the invariants in \rr{nsym-1a} and \rr{nsym-1b}.

Hence, the irreducible basis for $\tH$ contains, at most,  $9$ tensors,  $\tv_i\ot\tv_j$.
In the case when $\tH$ is symmetric, the irreducible basis contains at most $6$ symmetric tensors given in \rr{gv-sm3}.
In ths case when $\tH$ is a skew-symmetric tensor , the irreducible basis contains at most $3$ skew-symmetric tensors, i.e.
\be \tv_i\ot\tv_j - \tv_j\ot\tv_i \, \spcc( i=1,2; j=2,3, i<j) \, . \ee

Alternatively, for $M\ge 1$ and $N \ge 0$, using the singular value decomposition
\be \tH_1 = \sum_{i=1}^3 \ld_i \tv_i\ot\tu_i \, , \ee
 we can easily prove that
\be \tH(\tH_t, \tA_r,\ta_s) = \sum_{i,j=1}^3 \hat{h}_{ij} \tv_i\ot\tu_j \, , \ee
where $\tu_j$ are the unit eigenvectors of $\tH_1^T\tH_1$, $\tv_i$ are the unit eigenvectors of $\tH_1\tH_1^T$ and the invariants
\be \hat{h}_{ij} = \tv_i \cdot \tH \tu_j \ne \hat{h}_{ji} \, \ee
are functions of the $9M+6N+3P-3$ invariants
\be \ld_i \, , \spcc \tu_i\cdot\tv_i \, , \spcc \tv_i\cdot\tH_t\tu_j \spcc (t\ge 2) \, , \spcc \tv_i\cdot\tA_r\tu_j \, , \spcc
\ta_s \cdot\tv_i \, . \ee

Smith \cite{smith71} claimed for a symmetric tensor $\tH$ and skew-symmetric tensors $\tH_t=\tW_t$, the irreducible basis for symmetric $\tH$ contains the set of symmetric tensors
\[ \tI \, , \spcc \tA_i \, , \spcc \tA^2_i \, , \spcc \tA_i\tA_j+ \tA_j\tA_i \, , \spcc
\tA_i^2\tA_j+ \tA_j\tA_i^2 \, , \spcc \tA_i\tA_j^2+ \tA_j^2\tA_i \]
\[ \ta_m\ot\ta_m \, , \spcc \ta_m\ot\ta_n +\ta_n\ot\ta_m \, , \spcc
\ta_m\ot\tA_i\ta_m + \tA_i\ta_m\ot\ta_m \, , \spcc  \ta_m\ot\tA_i^2\ta_m + \tA_i^2\ta_m\ot\ta_m \, , \]
\[
\tA_i(\ta_m\ot\ta_n - \ta_n\ot\ta_m) - (\ta_m\ot\ta_n - \ta_n\ot\ta_m)\tA_i \, , \]
\[ \tW_p^2 \, , \spcc \tW_p\tW_q+\tW_q\tW_p \, , \spcc \tW_p\tW_q^2-\tW_q^2\tW_p \, , \spcc \tW_p^2\tW_q-\tW_q\tW_p^2 \, , \] 
\[  \tA_i\tW_p - \tW_p\tA_i \, , \spcc \tW_p\tA_i\tW_p \, , \spcc  \tA_i^2\tW_p - \tW_p\tA_i^2 \, , \spcc
\tW_p\tA_i\tW_p^2 - \tW_p^2\tA_i\tW_p \, , \]
\[ \tW_p\ta_m \ot\tW_p\ta_m \, , \spcc \ta_m\ot \tW_p\ta_m  + \tW_p\ta_m\ot \ta_m \, , \spcc 
\tW_p\ta_m \ot \tW_p^2\ta_m  + \tW_p^2\ta_m \ot \tW_p\ta_m  \, , \]
\be\label{nsym-5} 
\tW_p(\ta_m\ot\ta_n - \ta_n\ot\ta_m) + (\ta_m\ot\ta_n - \ta_n\ot\ta_m)\tW_p \, , \ee
where $(i,j=1,2, \ldots , N; i<j)$, $(p,q = 1,2, \ldots , M ; p < q)$ and $(m,n =1,2, \ldots , P; m< n)$. 
In \rr{nsym-5}, it is clear that there is a large number of "complicated" symmetric tensors in the Smith \cite{smith71} irreducible basis and this
number is far greater than $6$, the number of symmetric tensors in our irreducible basis. We note that all of Smith's symmetric tensors in
\rr{nsym-5} can be expressed in terms of the six symmetric tensors given in \rr{gv-sm3}.

\section{Potential Vectors and Tensors}\label{sec-der-1}
In this Section, we consider vectors and tensors that can be obtained from differentiating a scalar-valued isotropic function
$W$, i.e.,
\be \tg = \pdf{W}{\ta} \, , \spcc \tG = \pdf{W}{\tV} \, , \spcc \tH = \pdf{W}{\tF} \, , \ee
where $\ta$ is a vector, $\tV$ is a symmetric tensor and $\tF$ is a non-symmetric tensor. We called these vectors/tensors, potential vectors/tensors.
For example, in non-linear
hyper-elasticity, the potential nominal stress ${\dis \tS=\pdf{\rp{W}{e}}{\tF} }$, where
$\tF$ is the deformation gradient tensor and $\rp{W}{e}$ is the strain energy function.
\subsection{Vector}
Let $W(\tH_t, \tA_r,\ta_s)$,  be a scalar-valued isotropic function and let
$\ta=\ta_1$. From Appendix B and following the work of Shariff \cite{shariff17a}, we obtain the relation
\[ \tg(\tH_t, \tA_r,\ta_s) = \pdf{W}{\ta} = \pdf{W}{\ld}\tv_1 + \left(\frd{1}{\ld} \pdf{W}{\tv_1} \cdot\tv_2 \right)\tv_2 + 
\left( \frd{1}{\ld} \pdf{W}{\tv_1} \cdot\tv_3\right) \tv_3  \]
\be\label{vect-qq1} = \pdf{W}{\ld}\tv_1  +  
\frd{1}{\ld} \left[ (\tI - \tv_1\ot\tv_1)^T \pdf{W}{\tv_1} \right] \, , \ee
where $\ld = \sqrt{\ta\cdot\ta}$. It is clear from \rr{vect-qq1}, since the coefficients of $\tv_i$ are scalar-valued isotropic functions,
$\tg$ is a vector-valued isotropic function.

\subsection{Symmetric Tensor-Valued Isotropic Function $\tG$}
In this case, we let ${\dis \tV=\tA_1 = \sum_{i=1}^3 \ld_i \tv_i\ot\tv_i}$. Shariff \cite{shariff17a}  has shown that tensor-valued isotropic function
\[ \tG(\tH_t, \tA_r,\ta_s) = \pdf{W}{\tV} \] 
\be\label{sym-da} =\sum_{i=1}^3 \pdf{W}{\ld_i} \tv_i\ot\tv_i  + \sum_{i,j=1 \, ,  i<j }^3 \frd{1}{2(\ld_i-\ld_j)}
(\pdf{W}{\tv_i}\cdot\tv_j - \pdf{W}{\tv_j}\cdot\tv_i)(\tv_i\ot\tv_j + \tv_j\ot\tv_i) \, . \ee

\subsection{Non-symmteric Tensor-Valued Isotropic Function  $\tH$}
In this case, in view of singular value decomposition, we have  ${\dis \tH_1=\tF =\sum_{i=1}^3 \ld_i \tv_i\ot\tu_i }$, where
$\ld_i$ are the square root of the eigenvalues of $\tF\tF^T$, $\tv_i$ is a unit eigenvector of $\tF\tF^T$ and
$\tu_i$ is a unit eigenvector of $\tF^T\tF$. Shariff \cite{shariff17a} (using a derivative convention used in Itskov \cite{itskov2013}) has shown that
tensor-valued isotropic function 
\[ \tH(\tH_t, \tA_r,\ta_s) = \pdf{W}{\tF} \]
\be = \sum_{i=1}^3 \pdf{W}{\ld_i}\tv_i\ot\tu_i + \sum_{i,j=1, i\ne j}^3 
\frd{\left(\ld_i (\pdf{W}{\tu_i}\cdot\tu_j - \pdf{W}{\tu_j}\cdot\tu_i) + \ld_j(\pdf{W}{\tv_i}\cdot\tv_j - \pdf{W}{\tv_j}\cdot\tv_i)
\right)
\tv_i\ot\tu_j}{\ld_i^2-\ld_j^2} \, .\ee

\section{Remark}
In this communication we have shown that we need only $3$ linearly independent vectors to represent both potential and non-potential vectors and 
and a maximum of only $9$ linearly independent tensors to represent both potenstial and non-potenstial tensors.
However, the number of functions in a Smith \cite{smith71} or Boehler \cite{boehler77}  irreducible basis required to represent a potential vector/tensor 
is generally not the same as that required to represent a non-potential vector/tensor. For example, consider finite strain transversely isotropic elasticity with the preferred direction $\ta$ in the undeformed
configuration. Let $\tS(\tC,\tL)$ be the second Piola-Kirchhoff stress tensor, where $\tC$ is the right Cauchy-Green tensor and
$\tL=\ta\ot\ta$. Using Smith \cite{smith71} and Boehler \cite{boehler77}  tensor functions, we have 
\be\label{rem-Sa} \tS = \alpha_0 \tI + \alpha_1\tL + \alpha_2 \tC + \alpha_3\tC^2 + \alpha_4(\tC\tL + \tL\tC) + \alpha_5(\tC^2\tL + \tL\tC^2) \, , \ee
where $\alpha_0 - \alpha_5$ are isotropic invariants of the set $\{\tC,\tL\}$.
For an hyperelastic material, there exist a strain energy function 
\be W(\tC,\tL) = \hat{W}(I_1,I_2,I_3,I_4,I_5) \, , \ee
where the invariants
\be I_1 =\tr \tC \, , \spcc I_2=\tr \tC^2 \, , \spcc I_3=\tC^3  \, , \spcc I_4=\tr (\tC\tL) \, , \spcc I_5=\tr (\tC^2\tL) \, . \ee 
The second (potential) Piola-Kirchhoff stress tensor then has the relation
\be\label{rem-S} \tS = \pdf{W}{\tE} = 2\pdf{\hat{W}}{\tI_1} \tI + 4\pdf{\hat{W}}{\tI_2}\tC + 6\pdf{\hat{W}}{\tI_3}\tC^2+
2\pdf{\hat{W}}{\tI_4}\tL + \pdf{\hat{W}}{\tI_5}(\tC\tL + \tL\tC) \, , \spcc \tE=\frd{1}{2}(\tC-I) \, . \ee
Comparing \rr{rem-Sa} and \rr{rem-S}, we observe that the representation for the hyperelastic material does not include the last term
in \rr{rem-Sa}, i.e., $\tC^2\tL + \tL\tC^2$. It seems {\it on the onset}, if we use Smith \cite{smith71} and Boehler \cite{boehler77}  irreducible functions,
 the constitutive equation
 \rr{rem-Sa} {\it cannot} be described by a strain energy function (see comments made in Itskov \cite{itskov2013} page 144). 
However, if we express the tensors 
\be  \tI \, , \spcc \tL \, , \spcc  \tC \, , \spcc  \tC^2 \, , \spcc \tC\tL + \tL\tC \, , \spcc \tC^2\tL + \tL\tC^2 \,  \ee
in terms of the tensors $\tv_i\ot\tv_j$ ($\tv_i$ is an eigenvector of $\tC$), where their scalar coefficients are isotropic invariants of
the set $\{\tC,\tL \}$, we could easily
equate \rr{rem-Sa} with \rr{rem-S}; which suggest that, when express in terms of the basis functions $\tv_i\ot\tv_j$,  the constitutive equation
 \rr{rem-Sa} {\it can} be described by a strain energy function. 

In general, following the above example, it can be easily shown that a non-potential vector/tensor can always be represented by a 
potential vector/tensor.

\section*{Appendix A: $P$-property}
\def\theequation{A\arabic{equation}}
\setcounter{equation}{0}

The description of the $P$-property uses the eigenvalues $\ld_i$ and eigenvectors $\tv_i$ of the symmetric tensor $\tA_1$ . A general anisotropic 
invariant, where its arguments are expressed in terms spectral invariants with respect to the basis
$\{ \tv_1,\tv_2, \tv_3 \}$ can be written in the form
\be\label{pa1}\Phi &=& \bar{W}(\ld_i, \tv_i\cdot\tA_r\tv_j, \tv_i\cdot\ta_s) \nn \\
&=& \tilde{W}(\lambda_1,\lambda_2,\lambda_3,\tv_1,\tv_2,\tv_3) \, , \ee
where
\be r=2,\ldots, M, \, , \spcc s=1,2,\ldots, P, \ee
and, in Eqn. \rr{pa1}$_2$, the appearance of $\tA_r$ and $\ta_s$ is suppressed to facilitate the description of the $P$-property.
$\tilde{W}$ must satisfy the symmetrical property
\be\label{pa1a}
\tilde{W}(\lambda_1,\lambda_2,\lambda_3,\tv_1,\tv_2,\tv_3) = \tilde{W}(\lambda_2,\lambda_1,\lambda_3,\tv_2,\tv_1,\tv_3) =
 \tilde{W}(\lambda_3,\lambda_2,\lambda_1,\tv_3,\tv_2,\tv_1) \, . \ee
In view of the non-unique values of $\tv_i$ and $\tv_j$ when $\lambda_i=\lambda_j$, a function $\tilde{W}$ should be independent of $\tv_i$ and $\tv_j$ when $\lambda_i=\lambda_j$,
and $\tilde{W}$ should be independent of $\tv_1$, $\tv_2$ and $\tv_3$ when $\lambda_1=\lambda_2=\lambda_3$. Hence, when  two or three of the principal stretches have equal values
the scalar function $\Phi$ must have any of the following forms
\[ \Phi = 
\left\{ \begin{array}{cc}  \rp{W}{a}(\ld,\ld_k,\tv_k) \, , &  \ld_i=\ld_j=\ld \, , i\ne j \ne k \ne i \\
\rp{W}{b}(\ld) \, , &  \ld_1=\ld_2=\ld_3=\ld \hspace*{\fill} \end{array}\right. \] 
For example, consider 
\be\label{apa-11} \Phi=\ta\tA_1\ta =\sum_{i=1}^3 \ld_i(\ta\bl\tv_i)^2  \, , \ee
where $\ta$ is a fixed unit vector and 
\be  \sum_{i=1}^3 (\ta\bl\tv_i)^2 = 1 \, . \ee . 
If
\be \ld_1=\ld_2=\ld \, , \ee
 we have 
\be  \Phi= \rp{W}{a}(\ld,\ld_3,\tv_3) =\ld + (\ld_3 - \ld)(\ta\bl\tv_3)^2 \, \ee
and in the case of $\ld_1=\ld_2=\ld_3=\ld$, 
\be \Phi=\rp{W}{b}(\ld) = \ld \, . \ee
Hence, the invariant \rr{apa-11} satisfies the $P$-property and we note that all the classical invariants described in Spencer \cite{spencer71} 
 satisfy the $P$-property.
In reference \cite{shariff20}, the  $P$-property described  here is extended to non-symmetric tensors such as the two-point deformation tensor $\tF$.

\section*{Appendix B}
\def\theequation{B\arabic{equation}}
\setcounter{equation}{0}
A dyadic product $\ta\ot\ta$ has the spectral representation
\be\label{apb-1} \ta\ot\ta = \ld \tv_1 \, , \spcc \ld = \sqrt{\ta\cdot\ta} \, , \spcc \tv_1 = \frd{1}{\ld}\ta \, .  \ee
The unit eigenvectors $\tv_2$ and $\tv_3$, associated with zero eigenvalues, are non-unique. In view of \rr{apb-1}, we have
\be\label{apb-2} d\ta = d\ld \tv_1 + \ld d\tv_1  = d\ld \tv_1 + \ld ( da_2 \tv_2 + da_3\tv_3) \, . \ee
Note that the above expression, have used the relation, for arbitrary,
\be\label{apb-4}
d\tv_1 = da_2 \tv_2 + da_3\tv_3 \, , \ee
where $da_2$ and $da_3$  are arbitrary. We can write
\be\label{apb-3} d\ta = \sum_{i=1}^3 (d\ta)_i \tv_i \, , \spcc (d\ta)_1 = d\ld \, , \spcc (d\ta)_2= \ld da_2 \, , \spcc (d\ta)_3=\ld da_3 \, . \ee
For a scalar isotropic function $W=\rp{W}{a}(\ta) =\rp{W}{s}(\ld,\tv_1)$. 
Express
\be\label{apb-5} \pdf{\rp{W}{a}}{\ta} = \sum_{i=1}^3 \left(  \pdf{\rp{W}{a}}{\ta} \right)_i \tv_i \, , \spcc 
\left(  \pdf{\rp{W}{a}}{\ta} \right)_i = \pdf{\rp{W}{a}}{\ta}\cdot\tv_i \, . \ee
We then have
\be dW = \sum_{i=1}^3 \left(  \pdf{\rp{W}{a}}{\ta} \right)_i (d\ta)_i  = 
\pdf{\rp{W}{s}}{\ld} d\ld + \pdf{\rp{W}{s}}{\tv_1}\cdot d\tv_1 \, . \ee
Using \rr{apb-3} to \rr{apb-5} and since $d\ld , da_2$ and $da_3$ are arbitrary, we obtain the relations
\be \left(  \pdf{\rp{W}{a}}{\ta} \right)_1 = \pdf{\rp{W}{s}}{\ld} \, , \spcc
\left(  \pdf{\rp{W}{a}}{\ta} \right)_2 = \frd{1}{\ld} \pdf{\rp{W}{s}}{\tv_1} \cdot\tv_2 \, , \spcc
\left(  \pdf{\rp{W}{a}}{\ta} \right)_3 = \frd{1}{\ld} \pdf{\rp{W}{s}}{\tv_1} \cdot\tv_3 \, . \ee

\end{document}